\documentclass{amsart}
\usepackage{amsmath,amssymb,amsthm,xspace}
\newcommand{\tiff}{if\textcompwordmark f\xspace}
\newcommand{\N}{\mathbb{N}}
\newcommand{\Z}{\mathbb{Z}}
\newcommand{\R}{\mathbb{R}}
\newcommand{\inv}{^{-1}}
\newcommand{\orth}{^\perp}
\newcommand{\ab}{_{ab}}
\newcommand{\set}[1]{\left\{#1\right\}}
\newcommand{\setst}[2]{\set{\,#1\mid#2\,}}
\renewcommand{\span}[1]{\left<#1\right>}
\newcommand{\spanst}[2]{\span{\,#1\mid#2\,}}
\newcommand{\abs}[1]{\left|#1\right|}
\newcommand{\Abs}[2][]{{\left\|#2\right\|_{#1}}}
\newcommand{\gint}[1]{\left\lfloor#1\right\rfloor}
\newcommand{\onto}{\twoheadrightarrow}
\newcommand{\func}[4][\to]{#2\colon#3#1#4}
\newtheorem{prop}{Proposition}
\newtheorem{lem}[prop]{Lemma}

\theoremstyle{definition}
\newtheorem{defin}{Definition}
\theoremstyle{remark}
\newtheorem*{rk}{Remark}
\begin{document}
\title[Deep pockets for all generating sets]{A group with deep pockets
for all finite generating sets}
\date{\today}
\author[A.D. Warshall]{Andrew~D. Warshall}
\address{Yale University\\Department of Mathematics\\P.O. Box
208283\\New Haven, CT 06520-8283\\USA}
\email{andrew.warshall@yale.edu}
\thanks{We thank our advisor, Andrew Casson, for his helpful
comments.}
\subjclass[2000]{20F65}
\begin{abstract}
We show that the discrete Heisenberg group has unbounded dead-end
depth with respect to every finite generating set. We also show that,
in contrast, it has bounded retreat depth.
\end{abstract}
\maketitle

\section{Introduction}
Let $G$ be any group and let $A$ be a generating set for $G$. Let
$\Abs[A]{g}$ denote the length of the minimal-length word in $A$
representing $g$. By the \emph{depth} (or \emph{dead-end depth}) of an
element $g\in G$ with respect to $A$ we mean the distance (in the word
metric with respect to $A$) from $g$ to the nearest $g'\in G$ with
$\Abs[A]{g'}>\Abs[A]{g}$. If there is no such $g'$, then the depth of
$g$ is infinite.

If $G$ is a finite group then it will have elements of infinite
depth. In contrast, this cannot happen for $G$ infinite, provided $A$
is finite. It is natural then to ask whether depth is bounded over all
group elements. The first result in this direction was due to
Bogopol'sk\u{i}i, who showed in \cite{Bo} that every infinite
hyperbolic group has a bound on the depth of its elements with respect
to any given finite generating set. Later, we showed the same in
\cite{W} for all infinite Euclidean groups, and Lehnert then showed it
in \cite{L} for all groups with more than one end.

It is not difficult, however, to construct infinite finitely generated
groups $G$ and generating sets $A$ such that there is no bound on the
depth of elements of $G$ with respect to $A$. The first example, given
by Cleary and Taback in \cite{CT}, was the lamplighter group, namely
$\Z_2\wr\Z=\spanst{a,t}{t^2,[a,a^{t^i}],i\in\N}$, with respect to the
generating set $\set{a,t}$. Although this group is not finitely
presented, this is not essential; a (more complicated) finitely
presented example was given by Cleary and Riley in
\cite{CR}\footnote{The published version of this paper contained an
error, pointed out by Lehnert. A corrected version is available from
the \texttt{arXiv}.}.

However, the property of having unbounded depth (also known as the
\emph{deep pockets} property) is not a generating-set invariant; this
was shown in joint work between this author and Riley in
\cite{RW}\footnote{The published version of this paper contained an
analogous error, which became clear after the error in \cite{CR} was
pointed out; however, the proof of Theorem~3 contained in Sections~2
and~3 is unaffected.}. Thus the question remains of whether there
exists an infinte group which has unbounded depth with respect to
every finite generating set. In this paper, we answer this question in
the affirmative.

In particular, let $H$ denote the discrete Heisenberg group
\[
\spanst{x,y}{[x,[x,y]],[y,[x,y]]}.
\]

Then we will prove the following

\begin{prop}\label{infinitedepth}
Let $A$ be a finite generating set for $H$. Then there is no bound on
the depth of all elements of $H$ with respect to $A$.
\end{prop}

\begin{defin}
Let $G$ be a group, $A$ a generating set for $G$. Let $g\in G$ be such
that $\Abs[A]{g}=l$. Then the \emph{retreat depth} of $g$ with respect
to $A$ is the minimal $d$ such that $g$ lies in an unbounded component
of $B_{l-d}(1)$. (Lehnert in \cite{L} refers to a concept similar to
this as \emph{strong depth}.)
\end{defin}

The above definition is motivated by Bowditch's Question~8.4 in
\cite{Be}, which asks in effect whether, for every infinite $G$ and
finite $A$, retreat depth is bounded for all $g\in G$. (This question
was answered in the negative by Erschler, using the same example of
the lamplighter group.) It is clear that if (with respect to some
fixed generating set) the depth of elements of a group is bounded,
then so is the retreat depth. The converse, however, fails fairly
resoundingly for the Heisenberg group.

\begin{prop}\label{retbound}
Let $A$ be a finite generating set for $H$. Then there exists $r\in\N$
such that, for every $h\in H$, the retreat depth for $h$ is $\le r$.
\end{prop}

The paper is organized is follows. In the brief Section~\ref{prel}, we
will remind the reader of some facts about $\Z^2$ and metric geometry,
the proofs of which can be found in \cite{W}, and fix some notations.
In Section~\ref{id}, we will give the proof of
Proposition~\ref{infinitedepth}. Finally, in Section~\ref{rb}, we will
give the proof of Proposition~\ref{retbound}.

\section{Preliminaries}\label{prel}
To prove either of the main results of this paper, we will first need
to recall some facts about the geometry of $\Z^2$ with respect to
arbitrary finite generating sets. Let $A$ be such a generating set and
let $B$ be the (closed) convex hull of $A\cup A\inv$. Then $B$ is a
polygon in $\R^2$ and we have the following

\begin{prop}\label{geo}
Let $a_1$, $a_2\in A\cup A\inv$ be adjacent vertices of $B$ and let
$i_1$, $i_2\in\Z_{\ge 0}$. Then $i_1a_1+i_2a_2$ is a geodesic word in
$A$. In particular, there are $C$ and $D\in\N$ depending only on $A$
such that every $v\in\Z^2$ is within $C$ in the standard $L^1$ norm
and $D$ in the norm with respect to $A$ of some $v'\in\Z^2$
represented by a geodesic word of this form. It follows that
$\Abs[A]{v}-D\le\Abs[A]{v'}\le\Abs[A]{v}+D$ and that, if
$\Abs[B]{\cdot}$ is taken to mean the norm whose unit ball is $B$,
then $\abs{\Abs[A]{v}-\Abs[B]{v}}$ is bounded independently of $v$.
\end{prop}

We omit the proof; for details see \cite{W}.

For the first result we will also need an auxiliary fact, having
nothing to do with groups.

\begin{prop}\label{fuzz}
Let $f$ be a function from a metric space $A$ to $\Z$ and $n\in\Z$.
Suppose there exists $a\in A$ and $r\in\Z$ such that for all $a'\in
B_r(a)$, $f(a')\le f(a)+n$. Then there exists some $a'\in A$ such that
$f$ attains a maximum on $B_{r/n}(a')$ at $a'$.
\end{prop}

The proof of this proposition may also be found in \cite{W}. Its
relevance is that, to prove that the depth of a group is infinite, we
may replace distance from the identity with any other function
differing from it by at most a finite additive constant.

We agree to use $x$ and $y$ to denote the two standard generators of
$H$ and $\func[\onto]{\phi}{H}{H\ab}$ to denote the abelianization
map.

\section{Proof of Proposition~\ref{infinitedepth}}\label{id}
We begin by showing the existence of a particularly nice collection of
words of ``almost minimal'' length.

\begin{prop}\label{cobound}
Let $A$ be a generating set for $H$. Then there exist $G$, $I\in\R$,
$i_1$, $i_2$, \dots $\in\N\cup\set{0}$ and a sequence of words
$w_{00}$, $w_{10}$, \dots, $w_{1i_1}$, $w_{20}$, \dots, $w_{2i_2}$,
$w_{30}$, \dots in the letters of $A\cup A\inv$ such that
\begin{itemize}
\item $l(w_{ni})=2n$ for all $n$ and $i$,
\item for all $n$, $w_{(n+1)0}$ is obtained from $w_{ni_n}$ by
inserting a letter of $A$ and its inverse at positions in $w_{ni_n}$
separated by one letter,
\item for all $n>0$ and $i$, $w_{n(i+1)}$ is obtained from $w_{ni}$ by
interchanging one of the letters inserted to produce $w_{n1}$ with an
adjacent letter,
\item each $w_{ni}$ represents $[x,y]^{k_{ni}}$, where for every
$m\in\N$ there are $n$ and $i$ such that $\abs{m-k_{ni}}\le G$ and
\item for all $n$ and $i$, $l(w_{ni})-\Abs[A]{[x,y]^{k_{ni}}}\le I$.
\end{itemize}
\end{prop}

For $A$ a generating set for $H$ and $n\in\Z_{\ge0}$, let
\[
I_A(n)=\max\setst{\abs{k}}{\Abs[A]{[x,y]^k}=n}.
\]

Proposition~\ref{cobound} will now follow from

\begin{prop}\label{fattest}
Let $A$ be a generating set for $H$. Then there exist $D$, $E$ and
$F\in\R$ with $E>0$ and for each $n\in\N\cup\set{0}$ a word $w_n$ in
the letters of $A$ and their inverses such that
\begin{itemize}
\item $l(w_n)=2n$,
\item each $w_n$ is obtained from $w_{n-1}$ by inserting a letter of
$A$ and its inverse (not in general adjacently to each other) and
\item each $w_n$ represents $[x,y]^{k_n}$, where $I_A(2n)-Dn\le k_n\le
I_A(2n)$ and $En^2-Fn\le k_n\le En^2$.
\end{itemize}
\end{prop}

We postpone the proof.

\begin{proof}[Proof of Proposition~\ref{cobound} assuming Proposition~\ref{fattest}]
We proceed inductively on the assumption that
\[
w_{(n-1)i_{n-1}}=w_{n-1}
\]
from Proposition~\ref{fattest}. Let the two letters added to $w_{n-1}$
to make $w_n$ be $a$ and $a\inv$. Let $w_{n1}$ be obtained by
inserting $a$ and $a\inv$ at positions separated by one letter. Then,
let $w_{n2}$, \dots, $w_{ni_n}$ be obtained by interchanging either
$a$ or $a\inv$ with an adjacent letter at each step to yield
eventually $w_{ni_n}=w_n$. Note that we can arrange that
$i_n\le2n$. This construction clearly fulfills the first three
conditions.

Since each $w_{ni}$ contains each letter and its inverse an equal
number of times it represents some power of $[x,y]$. But each step in
the above construction changes the exponent by at most
\[
M'_A=\max\setst{\abs{k}}{[a_1,a_2]=[x,y]^k,a_1,a_2\in A\cup A\inv}.
\]
Since $w_{00}$ represents the identity and $w_{ni_n}=w_n$ represents
$[x,y]^{k_n}$ where $k_n\ge En^2-Fn$ and $E$ and $F$ are independent
of $n$ as in Proposition~\ref{fattest}, the $k_{ni}$ must increase
without bound. Thus the fourth condition is proven if we let $G=M'_A$.

Set $k_{max}(A)=\max\setst{\abs{k}}{x^iy^j[x,y]^k\in A}$. Each
$w_{ni}$ ($n>0$) represents
\[
[x,y]^{k_{ni}},
\]
where
\begin{multline*}
k_{ni}\ge k_{n-1}-2nM'_A\ge
I_A(2(n-1))-Dn-2nM'_A\\\ge4(n-C-1)^2M_A-2(n+C-1)k_{max}(A)-Dn-2nM'_A.
\end{multline*}
For $n$ sufficiently large, letting
$J=(D+2M'_A+2k_{max}(A))/(4M_A)+C+2$, this is
$\ge4(n-J)^2M_A+2(n-J)k_{max}(A)\ge I_A(2n-2J)$, where the second
inequality holds by the definitions of $I_A$, $k_{max}(A)$ and
$M_A$. This proves the final condition, since the $w_{ni}$ are thus
within $2J$ of minimal length. (The finitely many insufficiently large
values of $n$ can only increase $I$ by a finite amount.)
\end{proof}

\begin{prop}\label{polygon}
Let $A$ be a generating set for $H$. Then there exist $C$ ($=G$ from
Proposition~\ref{cobound}), $D$ and $E\in\Z$ with the following
property. For every $n\in\Z$ there exists $n'\in\Z$ with
$\abs{n-n'}\le C$ such that, if $v=v_1\phi(x)+v_2\phi(y)\in H\ab$ with
$\Abs[\set{\phi(x),\phi(y)}]{v}\le\sqrt[6]{\abs{n}}$, then there exist
words $w$ and $w'$ in the letters of $A$ and their inverses such that
\begin{itemize}
\item $w$ represents $[x,y]^{n'}$
\item $l(w)\le\Abs[A]{[x,y]^{n'}}+D$,
\item $w$ can be obtained from some word in the sequence given by
Proposition~\ref{cobound} by cyclic permutation and possibly inverting
all letters,
\item $\phi(w')=v$,
\item $l(w)\ge l(w')-E$ and
\item $l(w'w\inv)\le\sqrt[5]{\abs{n'}}$.
\end{itemize}
\end{prop}

\begin{rk}
Let $w'$ represent $x^{v_1}y^{v_2}[x,y]^{n''}$. Then, since
$\Abs[A]{x^{v_1}y^{v_2}[x,y]^{n''-n'}}\le\sqrt[5]{\abs{n'}}$, it
follows that there is some $F$ depending only on $A$ such that
$\abs{n''-n}\le F\abs{n'}^{2/5}+C$.
\end{rk}

We will need the following proposition, which asserts basically that
minimal-length representatives of powers of $[x,y]$ cannot be too long
and thin.

\begin{prop}\label{dirvar}
Let $A$ be a generating set of $H$. Then there are $C$, $D$ and
$E\in\R$ such that the following is true.  Let $a_1\dots a_m$ be a
minimal-length word in the letters of $A$ such that $a_1\dots
a_m=[x,y]^n$. Let $\func{f}{H\ab}{\R}$ be linear of norm $1$.  Then
there is a subword $a_i\dots a_j$ of $a_1\dots a_m$ with
$C\sqrt{\abs{n}}\ge j-i\ge D\abs{\sqrt{n}}-E$ and
$C\sqrt{\abs{n}}\ge\abs{f(\phi(a_i\dots a_j))}\ge D\sqrt{\abs{n}}-E$.
\end{prop}

\begin{lem}\label{notstraight}
Let $A$ be a generating set of $H$. Then there are $C$, $D$ and
$E\in\R$ such that the following is true.  Let $a_1\dots a_m$ be a
minimal-length word in the letters of $A$ such that $a_1\dots
a_m=[x,y]^n$. Let $\func{f}{H\ab}{\R}$ be linear of norm $1$. Let
$P_n$ be $\set{e,a_1,a_1a_2,\dots,a_1\dots a_m=[x,y]^n}$. Then the
diameter of $f(\phi(P_n))$ is at most $C\sqrt{\abs{n}}$ and at least
$D\sqrt{\abs{n}}-E$.
\end{lem}

\begin{proof}
Let $M(A)=\max_{a\in A}\Abs[\set{x,y}]{a}$.  We know that there is
$F\in\R$ such that $m\le F\sqrt{\abs{n}}$, since this is so with
respect to the standard generating set. Thus we know that the diameter
of $P_n$ with respect to the standard generators of $H$ is at most
\[
mM(A)\le FM(A)\sqrt{\abs{n}}.
\]
Since neither $f$ nor $\phi$ increases distances (where we consider
$H\ab$ with respect to the images of the standard generators) the same
bound holds for the diameter of $f(\phi(P_n))$, so we may take
$C=FM(A)$.

Conversely, choose a minimal-length word in the standard generators
for each letter of $A$, so $a_1\dots a_m$ becomes $b_1\dots b_M$,
$b_i\in\set{x,y,x\inv.y\inv}$, $M\le M(A)m$. Then, if
\[
\func{g}{H\ab}{\R}=f\circ\begin{pmatrix}0&1\\-1&0\end{pmatrix},
\]
(where the matrix is with respect to the basis
$\span{\phi(x),\phi(y)}$) we have
\begin{multline*}
\abs{n}=\abs{\sum_{i=1}^Mg(b_i)[f(\phi(b_1\dots
b_{i-1}))+\frac{f(\phi(b_1))}{2}]}\\\le\sum_{i=1}^M\abs{g(b_i)}\cdot[\max_{i\in\set{1,\dots,M}}\abs{f(\phi(b_1\dots
b_{i-1}))}+\frac{\max_{b\in\set{x,y}}\abs{f(\phi(b))}}{2}\\\le
M[\max_{i\in\set{0,\dots,m}}\abs{f(\phi(a_1\dots
a_i))}+M(A)+\frac{1}{2}]\\\le
M(A)m[\max_{i\in\set{0,\dots,m}}\abs{f(\phi(a_1\dots a_i))}+2M(A)\\\le
FM(A)\sqrt{\abs{n}}[\max_{i\in\set{0,\dots,m}}\abs{f(\phi(a_1\dots
a_i))}+2M(A).
\end{multline*}
But the diameter of $f(\phi(P_n))$ is
$\ge\max_{i\in\set{0,\dots,m}}\abs{f(\phi(a_1\dots a_i))}$, so by the
above it is
\[
\ge\frac{\sqrt{\abs{n}}}{FM(A)}-2M(A).
\]
Thus we may take $D=1/(FM(A))$ and $E=2M(A)$.
\end{proof}

\begin{proof}[Proof of Proposition~\ref{dirvar}]
Lemma~\ref{notstraight} gives us $C$, $D$ and $E$ independent of $n$
and some subword $a_i\dots a_j$ such that
$C\sqrt{\abs{n}}\ge\abs{f(\phi(a_i\dots a_j))}\ge
D\sqrt{\abs{n}}-E$. But since $f$ and $\phi$ do not increase
distances, we know $j-i\ge D\sqrt{\abs{n}}-E-1$. Finally, $j-i<m\le
F\sqrt{\abs{n}}$ for some $F$ depending only on $A$. Since we can
absorb the $1$ into $E$ and redefine $C$ to be the greater of $C$ and
$F$, we are done.
\end{proof}

\begin{proof}[Proof of Proposition~\ref{polygon}]
Let $C=G$ and $D=I$ from Proposition~\ref{cobound}, and let $w$
representing $[x,y]^{n'}$ be the word given by that proposition with
$\abs{n-n'}\le C$. Then the first three conditions hold.

Let $v\in H\ab$ be such that
$\Abs[\set{\phi(x),\phi(y)}]{v}\le\sqrt[6]{\abs{n'}}$ but no $w'$
exists satisfying the remaining conditions. As in the discussion
before Proposition~\ref{geo}, identifying $H\ab$ with $\Z^2$, let $B$
be the closed convex hull of $\pm\phi(A)$ and let $\Abs[B]{\cdot}$ be
the norm with $B$ as its unit ball. Then Proposition~\ref{geo} gives
$F$ such that, for all $v'\in H\ab$,
$\abs{\Abs[B]{v'}-\Abs[\phi(A)]{v'}}\le F$. Express $w=w_1w_2$, where
$l(w_1)=\sqrt[5]{n'}/2-2F$. It follows (setting $E=2F$) that
$\Abs[B]{\phi(w_1)+v}=\Abs[B]{\phi(w_2)-v}>\sqrt[5]{n'}/2-F$. However,
everything we know about $w$ applies equally to any cyclic permutation
of $w$ and to any word obtained by inverting every letter of $w$.  We
conclude that, for any cyclic subword $u$ of $w$ with
$l(u)=\sqrt[5]{n'}/2-2F$, $\Abs[B]{\phi(u)\pm v}>\sqrt[5]{n'}/2-F$.

However, $\Abs[B]{\phi(u)}\le l(u)+F=\sqrt[5]{n'}/2-F$ for any such
subword $u$. Let $B'$ be the dilation of $B$ by
$\sqrt[5]{n'}/2-F$. Then $\phi(u)\in B'$ while $\phi(u)\pm v\notin
B'$. We refer to points on the boundary of $B'$ where the lower
(respectively upper) directional derivative with respect to $v$ of
$\Abs[B]{\cdot}$ is positive (resp.\ negative) as $v$-positive (resp.\
$v$-negative). Since, for any point in the interior of an edge, this
is equivalent to the outward normal to the edge's having positive
(resp.\ negative) dot product with $v$, we can speak of edges' being
$v$-positive or $v$-negative if we restrict attention to their
interior. Furthermore, a vertex of $B$ is $v$-positive (respectively
$v$-negative) \tiff both edges incident to it are.  In this language,
$\phi(u)$ lies within $\Abs{v}$ (in the Euclidean norm) of both a
point on a $v$-positive edge of $B'$ and one on a $v$-negative
edge.

Since $B'$ is convex and simple, there is $G$ such that $\phi(u)$ is
within $G\Abs{v}$ of a vertex of $B'$ at which a $v$-positive and a
$v$-negative edge meet. If there are any such vertices (as there must
bem by assumption), there are exactly two, and they lie opposite to
each other with respect to the origin. Call them $\pm c$ and let
$V\subset\Z^2$ be the subspace spanned by $c$. Let
$\func[\onto]{f}{\Z^2}{V\orth}$ be orthogonal projection, so that
$\ker f=V$. But then Proposition~\ref{dirvar} gives $I$ depending only
on $A$ such that for $n$ large enough there is a cyclic subword $u$ of
$w$ of length $\sqrt[5]{n'}/2-2F$ such that
$\Abs{f(\phi(u))}>I\sqrt[5]{n'}$. (If $n$ is not large enough, we just
set $n'=0$, increasing $C$ as needed.) But if we let
$t_{max}(A)=\max_{a\in A}\Abs{\phi(A)}$ then $\Abs{v}\le
t_{max}\sqrt[6]{n'}$, so $\Abs{f(\phi(u))}\le
G\Abs{v}+\Abs{f(c)}=G\Abs{v}\le Gt_{max}(A)\sqrt[6]{n'}$, which is a
contradiction for $n>(Gt_{max}(A)/I)^{30}+C$, so we are done if we let
$E=2F$.
\end{proof}

\begin{proof}[Proof of Proposition~\ref{infinitedepth}]
Let $n\in\N$ and let $n'$ be as in Proposition~\ref{polygon}. Let
$x^iy^j[x,y]^k$ be such that
$\Abs[A]{x^iy^j[x,y]^{k-n}}\le\sqrt[6]{n'}$. Let $v=(i,j)$ and let
$w$, $w'$, $n''$, $C$, $D$ and $E$ the words and constants given by
that proposition. Then
$l(w)\le\Abs[A]{[x,y]^{n'}}+D\le\Abs[A]{[x,y]^n}+C\Abs[A]{[x,y]}+D$.

By the remark following the statement of Proposition~\ref{polygon},
there is some $F$ depending only on $A$ such that
\[
\abs{n''-n}\le F\abs{n'}^{2/5}+C.
\]
Since $\Abs[A]{x^iy^j[x,y]^{k-n}}\le\sqrt[6]{n'}$, it follows that
$\Abs[\set{x,y}]{x^iy^j[x,y]^{k-n}}\le2\sqrt{G}\sqrt[6]{n'}$ for sone
$G>0$ depending only on $A$, whence $\abs{k-n}\le G\sqrt[3]{n'}$.
Thus 
\[
\abs{n''-k}\le\abs{n''-n}+\abs{n-k}\le
F\abs{n'}^{2/5}+G\sqrt[3]{n'}+C\le(F+G)\abs{n'}^{2/5}+C.
\]
Applying Proposition~\ref{cobound} to $k+n'-n''$ gives some $w''$
representing $[x,y]^{n'''}$ with $\abs{n'''+n''-n'-k}\le C$, so
\begin{multline*}
\abs{n'''-n}\le\abs{n'''+n''-n'-k}+\abs{n''-k}+\abs{n'-n}\\\le
C+(F+G)\abs{n'}^{2/5}+C+C=(F+G)\abs{n'}^{2/5}+3C.
\end{multline*}
It follows from Propositions~\ref{cobound} and~\ref{polygon} that a
cyclic permutation, combined with possibly inverting every element,
transforms $w''$ into a word that differs from $w$ only in changing
the position of boundedly many letters (say by $I$) and adding or
deleting at most $I$ letters. Assume without loss of generality that
$w''$ is so transformed. Decompose $w=w_1w_2$ and $w'$ as $w'_1w_2$,
so that $l(w_1)\le\sqrt[5]{n'}$ and $l(w'_1)\le l(w_1)+E$. Note that
$w'_1w_1\inv=w'w\inv=x^iy^j[x,y]^{n''-n'}$. Finally, decompose
$w''=w''_1w''_2$, where $l(w''_1)=l(w_1)$. (We may safely neglect the
case where $l(w'')<l(w_1)$, for then $n'-\sqrt[5]{n'}<I$, yielding a
bounded number of cases which may be ignored.) Let $t_{max}(A)$ again
denote $\max_{a\in A}\phi(a)$ and let $t'_{max}(A)$ denote
$\max(\Abs[A]{x},\Abs[A]{y})$.  Then $\Abs{\phi(w''_1)-\phi(w_1)}\le
It_{max}(A)$. Let
\[
M'_A=\max\setst{\abs{k}}{[a_1,a_2]=[x,y]^k,a_1,a_2\in A\cup A\inv}.
\]
Then there is a word $u'$ with $l(u')\le It_{max}(A)t'_{max}(A)$ with
$\phi(u')=\phi(w''_1)-\phi(w_1)$ and $w_1\inv w''_1=u'[x,y]^k$, where
$k\le l(u')l(w_1)M'_A$. This yields a word $u$ representing $w_1\inv
w''_1$ with
\begin{multline*}
l(u)\le It_{max}(A)t'_{max}(A)+8\sqrt{Gk}\\\le
It_{max}(A)t'{_max}(A)+8\sqrt{Gl(u')l(w_1)M'_A}\\\le
It_{max}(A)t'_{max}(A)+8\sqrt{GIt_{max}(A)t'_{max}(A)M'_A}\sqrt[10]{n'}.
\end{multline*}
Thus $w'_1uw''_2$ has length (for some $J$ depending only on $A$)
\begin{multline*}
\le l(w'_1)+l(u)+l(w''_2)\\\le
l(w_1)+E+It_{max}(A)t'_{max}(A)+8\sqrt{GIt_{max}(A)t'_{max}(A)M'_A}\sqrt[10]{n'}+l(w_2)+I\\=l(w)+8\sqrt{GIt_{max}(A)t'_{max}(A)M'_A}\sqrt[10]{n'}+It_{max}(A)t'_{max}(A)+E+I
\end{multline*}
and represents
$w'_1w_1\inv[x,y]^{n'''}=x^iy^j[x,y]^{n'''+n''-n'}$. Since
\[
\Abs{n'''+n''-n'-k}\le C,
\]
there is a word of length at most
\[
l(w)+8\sqrt{GIt_{max}(A)t'_{max}(A)M'_A}\sqrt[10]{n'}+It_{max}(A)t'_{max}(A)+E+I+C\Abs[A]{[x,y]}
\]
representing $x^iy^j[x,y]^k$. Since
\[
l(w)\le\Abs[A]{[x,y]^{n'}}+D\le\Abs[A]{[x,y]^n}+C\Abs[A]{[x,y]}+D,
\]
this length bound is $\le\Abs[A]{[x,y]^n}+K\sqrt[10]{n'}+L$, where
\[
K=8\sqrt{GIt_{max}(A)t'_{max}(A)M'_A}
\]
and $L=It_{max}(A)t'_{max}(A)+E+I+2C\Abs[A]{[x,y]}+D$ both depend only
on $A$. Since
$\lim_{n\to\infty}\sqrt[6]{n'}/(K\sqrt[10]{n'}+L)=\infty$, we are done
by Proposition~\ref{fuzz}.
\end{proof}

It remains to prove Proposition~\ref{fattest}. We begin with a
definition.

\begin{defin}
Let the \emph{isoperimetric constant} of a polygon in $\R^2$ with
respect to a given norm be the ratio of its enclosed area (in the
standard measure on $\R^2$) to the sum of its side lengths with
respect to the given norm (which we call the perimeter with respect to
that norm).
\end{defin}

\begin{lem}\label{areasym}
Let $\Abs{\cdot}$ be a norm on $\Z^2$ whose fundamental polygon is
actually a Euclidean polygon. Then there exists a polygon in $\R^2$
with maximal isoperimetric constant among all polygons in $\Z^2$ such
that
\begin{itemize}
\item it is convex and simple,
\item each side is parallel to another side of the same length and
\item every side is parallel to the ray from the origin to some vertex
of the fundamental polygon.
\end{itemize}
\end{lem}

\begin{proof}
Any polygon with maximal isoperimetric constant is clearly convex and
simple. For any such polygon, circumscribe it with a polygon each side
of which is parallel to the ray from the origin to some vertex of the
fundamental polygon. This polygon will have at most the same perimeter
and at least as large an area. Thus there exists a polygon with
maximal isoperimetric constant satisfying the first and third
conditions.

For the second condition, choose a vertex of the polygon and consider
the other point of the polygon which divides the perimeter into two
equal parts. A line segment connecting these two points must then
divide the area into two equal parts as well, for otherwise the
isoperimetric constant would not be maximal. Thus, if we replace one
half of the polygon with the other half rotated by $\pi$ and
translated, the area, hence the isoperimetric constant, will be the
same as before, so still maximal. The result must fulfill the third
condition if the original polygon did, and if it fails to be convex
and simple then it does not have maximal isoperimetric constant, so
neither did the original, a contradiction.
\end{proof}

\begin{lem}\label{intapprox}
Let $m\in\N$ and $b_1$, $b_2$, \dots, $b_m\in\R_+$ such that
$\sum_{i=1}^mb_i=1$. Then for every $n\in\N\cup\set{0}$ and
$i\in\set{1,\dots,m}$ there exists $b_{ni}\in\N\cup\set{0}$ such that
\begin{itemize}
\item $b_{0i}=0$ for all $i$,
\item $\sum_{i=1}^mb_{ni}=n$ for all $n$,
\item $b_{ni}$ is a nondecreasing function of $n$ for all $i$ and
\item $\abs{b_{ni}-nb_i}<m$ for all $n$ and $i$.
\end{itemize}
\end{lem}

\begin{rk}
It follows that for every $n$ there is some $j\in\set{1,\dots,m}$ such
that $b_{nj}=b_{(n-1)j}+1$ and $b_{ni}=b{(n-1)i}$ for all $i\ne j$.
\end{rk}

\begin{proof}
We construct the $b_{ni}$ inductively. First, let all the
$b_{0i}=0$. Assuming the $b_(n-1)i$ constructed, let
$j\in\set{1,\dots,m}$ be such that $nb_j-b_{(n-1)j}\ge
nb_i-b_{(n-1)i}$ for all $i\in\set{1,\dots,m}$. Then let
$b_{nj}=b_{(n-1)j}+1$ and $b_{ni}=b_{(n-1)i}$ for $i\ne j$. The first
three conditions are then clearly satisfied.

But let $n\in\N\cup\set{0}$ and define $j$ as above. We know
$\sum_{i=1}^m(nb_i-b_{(n-1)i})=1$, so $nb_j-b_{(n-1)j}>0$. Then
$nb_j-b_{nj}=nb_j-b_{(n-1)j}-1>-1$. But, for $i\ne j$,
$nb_i-b_{ni}=nb_i-b_{(n-1)i}\ge(n-1)b_i-b_{(n-1)i}>-1$ by induction,
since, for all $i$, $0b_i-b_{0i}=0>-1$. Since
$\sum_{i=1}^m(nb_i-b_{ni})=0$ for all $n$, this implies the fourth
condition.
\end{proof}

\begin{proof}[Proof of Proposition~\ref{fattest}]
If $A$ is a generating set for $H$, then $\phi(A)$ is a generating set
for $H\ab$, which we again identify with $\Z^2$. As in the discussion
before Proposition~\ref{geo}, let $B$ be the convex hull of
$\pm\phi(A)$ and let $\Abs[B]{\cdot}$ denote the norm with $B$ as its
unit ball.  Then Proposition~\ref{geo} gives $C\in\R$ such that
$\abs{\Abs[\phi(A)]{g}-\Abs[B]{g}}\le C$ for all $g\in\Z^2$.

Consider the polygon of maximal isoperimetric constant with respect to
$\Abs[B]{\cdot}$ given by Lemma~\ref{areasym}. Denote this maximal
isoperimetric constant by $M_A$. Choose $m$ sides of the polygon,
taken in order. Scale the polygon to have perimeter $2$. It must have
evenly many sides, so let $m$ be half its number of sides. Let $b_1$,
\dots, $b_m\in\R_+$ be the lengths (in the norm $\Abs[B]{\cdot}$) of
the chosen $m$ sides; the lengths of the other $m$ sides are the
same. Then $\sum_{i=1}^mb_i=1$, so choose $b_{ni}$ by
Lemma~\ref{intapprox}. Since each side of our original polygon was
parallel to the projection under $\phi$ of an element of $A$, we get
for each $n$ a word $w_n$ of length $2n$ (thus satisfying the first
condition)
\[
a_1^{b_{n1}}\dots a_m^{b_{nm}}a_1^{-b_{n1}}\dots
a_m^{-b_{nm}}=[x,y]^{k_n},
\]
where the $a_i\in A\cup A\inv$ are independent of $n$ and
$k_n\in\Z$. By the remark following Lemma~\ref{intapprox} the $w_n$
satisfy the second condition as well.

Each $k_n$ is (up to sign) the area of the polygon with pairs of
opposite sides parallel to those of the original polygon but with
lengths $b_{ni}$. This polygon will have perimeter $2n$. But then
\begin{multline*}
0\le4n^2M_A-\abs{k_n}\le\sum_{i=1}^m\abs{b_{ni}-nb_i}(n+\sum_{j=1}^m\abs{b_{nj}-nb_j})\\<\sum_{i=1}^mm(n+m^2)=m^2(n+m^2)
\end{multline*}

But if we set $k_{max}(A)=\max\setst{\abs{k}}{x^iy^j[x,y]^k\in A}$
then $I_A(2n)\le4(n+C)^2M_A+2(n+C)k_{max}(A)$. If $n=0$ then
$I_A(2n)=\abs{k_n}=n^2=0$, while if $n>0$ then
\begin{multline*}
0\le I_A(2n)-\abs{k_n}\\\le
m^2(n+m^2)+8nCM_A+4C^2M_A+2nk_{max}(A)\\\le(m^4+m^2+4C^2M_A+8CM_A+2k_{max}(A))n
\end{multline*}
and $m^2(n+m^2)\le(m^2+m^4)n$. We may thus take
$D=m^4+m^2+4C^2M_A+8CM_A+2k_{max}(A)$, $E=4M_A$ and $F=m^4+m^2$.
\end{proof}

This completes the proof of Proposition~\ref{infinitedepth}.

\section{Proof of Proposition~\ref{retbound}}\label{rb}
Proposition~\ref{retbound} will follow from

\begin{prop}\label{kout}
Let $A$ be a generating set for $H$. Then there exists $C\in\N$ with
the following property. Let $i$, $j$, $k_1$, $k_2\in\Z$ with either
$ij/2\le k_1\le k_2$ or $ij/2\ge k_1\ge k_2$. Then
\[
\Abs[A]{x^iy^j[x,y]^{k_1}}\le\Abs[A]{x^iy^j[x,y]^{k_2}}+C.
\]
\end{prop}

This result in turn follows from another result, the statement of
which requires the following

\begin{defin}
Let $A$ be a generating set for $H$. Let $i\phi(x)+j\phi(y)\in H\ab$
and $n$ be the minimal length with respect to $H$ of any element of
$\phi\inv(i\phi(x)+j\phi(y))$. For $n'\ge n$, let
$k^{max}_{(i,j)}(n')$ (respectively $k^{min}_{(i,j)}(n')$) be the
maximum (resp.\ minimum) $k$ such that there is a word of length $\le
n'$ representing $x^iy^j[x,y]^k$.
\end{defin}

\begin{prop}\label{allkapprox}
Let $A$ be a generating set for $H$. Then there are $D$ and $I\in\N$
with the following property. Let $i\phi(x)+j\phi(y)\in H\ab$, $n$ be
the minimal length with respect to $H$ of any element of
$\phi\inv(i\phi(x)+j\phi(y))$ and $n'\ge n+I$. Then for every $k$ with
$k^{min}_{(i,j)}(n')\le k\le k^{max}_{(i,j)}(n')$ there is $k'\in\Z$
with $\abs{k'-k}\le D$ such that $\Abs[A]{x^iy^j[x,y]^{k'}}\le
n'$. Furthermore we have $k^{min}_{(i,j)}(n')\le ij/2\le
k^{max}_{(i,j)}(n')$.
\end{prop}

We postpone the proof.

\begin{proof}[Proof of Proposition~\ref{kout}]
Let $D$ and $I\in\N$ be as in Proposition~\ref{allkapprox}. Let
$a_1\dots a_m$ be a minimal-length word in the letters of $A$ and
their inverses representing $x^iy^j[x,y]^{k_2}$. (Thus
$m=\Abs[A]{x^iy^j[x,y]^{k_2}}$.) It follows by definition that
\[
k_{(i,j)}^{min}(m+I)\le k_{(i,j)}^{min}(m)\le k_2\le
k_{(i,j)}^{max}(m)\le k_{(i,j)}^{max}(m+I)
\]
and that $m+I\ge n+I$.

Assuming without loss of generality that $ij/2\le k_1\le k_2$, we see
by the last sentence of Proposition~\ref{allkapprox} applied with
$n'=m+I$ that
\[
k_{(i,j)}^{min}(m+I)\le ij/2\le k_1\le k_2\le k_{(i,j)}^{max}(m+I).
\]
and thus that there exists a $k'\in\Z$ with $\abs{k'-k_1}\le D$ such
that
\[
\Abs[A]{x^iy^j[x,y]^{k'}}\le m+I.
\]
Thus, by the triangle inequality,
\[
\Abs[A]{x^iy^j[x,y]^{k_1}}\le\Abs[A]{x^iy^j[x,y]^{k'}}+\Abs[A]{[x,y]^D}\le
m+I+D\Abs[A]{[x,y]}.
\]
We are done if we let $C=I+D\Abs[A]{[x,y]}$.
\end{proof}

\begin{proof}[Proof of Proposition~\ref{retbound}]
For any $x^iy^j[x,y]^k\in H$, the sequence
\[
\setst{x^iy^j[x,y]^{k'}}{(k'-k)(k-ij/2)>0,\abs{k'-ij/2}>\abs{k-ij/2}}
\]
has distance from the identity (with respect to any finite generating
set) increasing without bound, since the metric with respect to any
finite generating set is proper. However, consecutive elements are at
$A$-distance $\Abs[A]{[x,y]}$ from each other, and, by
Proposition~\ref{kout} applied to $k_1=k$, $k_2=k'$, every element is
at distance at least $\Abs[A]{x^iy^j[x,y]^k}-C$ from the identity with
respect to $A$. Thus we may take $r=C$.
\end{proof}

We now proceed with the proof of Proposition~\ref{allkapprox}.  We will
first need two lemmas.

\begin{lem}\label{extremes}
Let $A$ be a generating set for $H$. Then there exists $C\in\N$ with
the following property. Let $n\in\N$ and $\set{a_1,\dots,a_n}$ be a
multiset of letters of $A$ and their inverses. Their projections to
$H\ab$ have a well-defined sum, which we denote
$i\phi(x)+j\phi(y)$. Then the average of the minimal and maximal $k$
such that some reordering $\sigma$ of the $a_i$ represents
$x^iy^j[x,y]^k$ is within $Cn$ of $ij/2$.
\end{lem}

\begin{proof}
Let $a_t=x^{b_t}y^{c_t}[x,y]^{k_t}$.  Averaging over all $n!$ possible
choices for $\sigma$ gives an average $k$-value of
\[
\frac{\sum_{\sigma\in S_n}\sum_{1\le s<t\le
n}c_{\sigma(s)}b_{\sigma(t)}}{n!}+\sum_{t=1}^nk_t=\frac{n!\sum_{1\le
s\ne t\le n}b_sc_t}{2n!}+\sum_{t=1}^nk_t,
\]
which differs (absolutely) from
$\sum_{t=1}^nb_t\cdot\sum_{t=1}^nc_t/2$ by
\[
\abs{\sum_{t=1}^n(k_t-b_tc_t)}\le n\max_{t=1}^n\abs{k_t-b_tc_t}\le
n\max_{x^by^c[x,y]^k\in A\cup A\inv}\abs{k-bc}.
\]
The average of the minimal and maximal $k$-values equals the average
$k$-value, since the distribution is symmetrically distributed about
the average, as may be seen by considering what happens to $k$ when
$\sigma$ is followed by the permutation which changes the order of
every pair of elements. The result follows, letting $C$ be an integer
greater than $\max_{x^by^c[x,y]^k\in A\cup A\inv}\abs{k-bc}$.
\end{proof}

\begin{lem}\label{kapprox}
Let $A$ be a generating set for $H$. Then there are $E\in\R_{>0}$ and
$F\in\N$ with the following property.  Let $i\phi(x)+j\phi(y)\in
H\ab$, $n$ be the minimal length with respect to $H$ of any element of
$\phi\inv(i\phi(x)+j\phi(y))$ and $d\ge0$.  Then there exists
$\set{a_1,\dots,a_{n+d}}$, a multiset of letters of $A$ and their
inverses, such that $\sum_{t=i}^{n+d}\phi(a_t)=(i,j)$ and the
difference between the maximal and minimal possible exponent of
$[x,y]$ obtainable by multiplying the $a_t$ in some order is $\ge
f_{E,F}(n,d)=\max(E\max(n-F,0)\max(d-F,0),Ed^2-Fd)$. In particular,
$k^{max}_{(i,j)}(n+d)-k_{min}(n+d)\ge f_{E,F}(n,d)$.
\end{lem}

\begin{proof}
We show that the difference can be made to exceed each of the two
expressions of which $f_{E,F}(n,d)$ is the maximum.

We will use $\Abs{\cdot}$ to denote the standard $L^1$ norm on $\Z^2$,
which we identify with $H\ab$.  It follows from Proposition~\ref{geo}
that there are $C$ and $D\in\N$ with the following property.  Let $w$
be a word of length $n$ representing an element $h\in\phi\inv((i,j))$
and $v=\phi(h)$. Then there is $v'$ with $\Abs{v'-v}\le C$ such that
there is an element $h'\in\phi\inv(v')$ with $n-D\le\Abs[A]{h'}\le
n+D$ and such that this length is attained by a word in $A\cup A\inv$
which is a product of letters projecting to at most two adjacent
vertices of the unit ball obtained from $A$ by that proposition. Let
$w'$ be this word representing $h'$. Since $\Z^2$ is not cyclic, the
unit ball must have at least four vertices, hence at least two pairs
of inverse vertices; we can thus choose such a pair (say $\pm p$,
$p=\phi(a)$, $a\in A$) such that at most half the letters of $w'$
project to $\pm p$. Write $w'=a_1\dots a_m$, where $n-D\le m\le n+D$
and the $a_t\in A\cup A\inv$. Let $k_t$ be such that
$[x,y]^{k_t}=[a,a_t]$.  The projections of the $a_t$, being at most
two adjacent vertices of the unit ball, must all lie in one of the two
closed half-planes bounded by $\span{p}$. Thus $k_t$ is either
nonnegative for every $t$ or nonpositive for every $t$. Assume without
loss of generality that it is nonnegative.  Furthermore, it is $0$ for
at most $m/2$ values of $i$ by our choice of $p$.  If $d>D$, then
$a^{\gint{(d-D)/2}}w'a^{-\gint{(d-D)/2}}$ represents
$h'[x,y]^{\gint{(d-D)/2}\sum_{t=1}^mk_t}$, while
$a^{\gint{(d-D)/2}}w'a^{\gint{(d-D)/2}}$ represents
$h'[x,y]^{-\gint{(d-D)/2}\sum_{t=1}^mk_t}$. Both these words have
length at most $n+D+2\gint{(d-D)/2}\le n+D+d-D=n+d$, and their
exponents of $[x,y]$ differ by $\gint{(d-D)/2}\sum_{t=1}^mk_t$.  If we
let
\[
c_{min}(A)=\min\setst{\abs{k}}{[a',a'']=[x,y]^k,a',a''\in A},
\]
we get that the two words' exponents of $[x,y]$ differ by
\begin{multline*}
\gint{\frac{d-D}{2}}\sum_{t=1}^mk_t\ge\frac{(d-D-1)c_{min}(A)m}{2}\\\ge
\frac{c_{min}(A)(d-D-1)(n-D)}{2}\ge\frac{c_{min}(A)(d-D-1)(n-D-1)}{2}.
\end{multline*}
Thus we are done with the first expression if we let $E=c_{min}(A)/2$
and $F=D+1$.

For the second expression, again let $w$, still of length $n$,
represent $x^iy^j[x,y]^k$. Let $w_x$ and $w_y$ be minimal-length words
in the letters of $A$ and their inverses representing $x$ and $y$,
respectively. Then, if we let $l_p=2(l(w_x)+l(w_y))$ we have
\[
w\left[w_x^{\gint{\frac{d}{l_p}}},w_y^{\gint{\frac{d}{l_p}}}\right]
\]
and
\[
w\left[w_y^{\gint{\frac{d}{l_p}}},w_x^{\gint{\frac{d}{l_p}}}\right]
\]
have length at most $n+d$ and represent
$x^iy^j[x,y]^{k+\gint{d/l_p}^2}$ and
$x^iy^j[x,y]^{k-\gint{d/l_p}^2}$, respectively.
But
\[
2\gint{\frac{d}{l_p}}^2\ge\frac{2d^2}{l_p^2}-\frac{4d}{l_p},
\]
so we are done with this expression if we let $E=2/l_p^2$ and
$F=4/l_p$.

Thus, we may simply take the lesser of these two values of $E$ and an
integer greater than or equal to these two values of $F$. The last
sentence then follows trivially.
\end{proof}

\begin{proof}[Proof of Proposition~\ref{allkapprox}]
For any word of length $n'$ in the $A$ and their inverses, we can
transform it to any permutation of itself by a succession of
transpositions of adjacent letters. Each such transposition changes
the exponent of $[x,y]$ by at most
\[
c_{max}(A)=\max\setst{\abs{k}}{[a_1,a_2]=[x,y]^k,a_1,a_2\in A\cup
A\inv}.
\]

Now consider a word $w_{max}$ of length at most $n'$ representing
\[
x^iy^j[x,y]^{k^{max}_{(i,j)}(n')}.
\]
By Lemma~\ref{extremes}, there is $C\in\N$ depending only on $A$ such
that we can permute the letters of $w_{max}$ so that it will represent
$x^iy^j[x,y]^{k_m}$ with $(k^{max}_{(i,j)}(n')+k_m)/2\le ij/2+Cn'$.
Thus $(k^{max}_{(i,j)}(n')+k^{min}_{(i,j)})/2\le ij/2+Cn'$.
Similarly, if $w_{min}$ is taken to represent
$x^iy^j[x,y]^{k^{min}_{(i,j)}(n')}$, then we can permute its letters
so it will represent some $x^iy^j[x,y]^{k_M}$ with
$(k^{max}_{(i,j)}+k_M)/2\ge ij/2+Cn'$, so
$(k^{max}_{(i,j)}(n')+k^{min}_{(i,j)})/2\ge ij/2-Cn'$. It follows that
$k^{max}_{(i,j)}-k_M\le4Cn'$ and $k_m-k^{min}_{(i,j)}\le4Cn'$. If (in
the notation of Lemma~\ref{kapprox}) $f_{E,F}(n,n'-n)\ge8Cn'$ for the
appropriate $E\in\R_{>0}$ and $F\in\N$ (which depend only on $A$)
then, by Lemma~\ref{kapprox},
$k^{max}_{(i,j)}(n')-k^{min}_{(i,j)}\ge8Cn'$, so the above
inequalities imply $k_M\ge k_m$. Thus there will exist $k'\in\Z$ with
$\abs{k'-k}\le c_{max}(A)/2$ such that some permutation of either
$w_{max}$ or $w_{min}$ represents $x^iy^j[x,y]^{k'}$, so we will be
done.

It remains to note that we will have $E(n-F)(n'-n-F)\ge8Cn'$ as soon
as $n'\ge E(n+F)(n-F)/[E(n-F)-8C]$. But, if $n\ge16C/E+2F$ we have
\begin{multline*}
\frac{E(n+F)(n-F)}{E(n-F)-8C}\le(n+F)\left[1+\frac{16C}{E(n-F)}\right]\\\le n+2F+\frac{16C}{E}\left(1+\frac{2F}{n}\right)\le n+2F+\frac{32C}{E},
\end{multline*}
where we use repeatedly that if $0\le x\le 1/2$ then
$0\le1/(1-x)\le1+2x$. If $n<16C/E+2F$, we will have
$E(n'-n)^2-F(n'-n)\ge8Cn'$ as soon as
$E(n'-n)^2-F(n'-n)\ge8C(n'-n+16C/E+2F)$. This is a quadratic
inequality in $n'-n$, which will hold so long as $n'-n\ge G$, say,
where $G$ depends only on $A$. So we are done if we let $D$ be an
integer $\ge c_{max}(A)/2$ and $I$ an integer $\ge\max(G,2F+32C/E)$.
The last sentence follows since if $f_{E,F}(n,n'-n)\ge8Cn'$ then it is
\emph{a fortiori} $\ge2Cn'$, so since $(k^{max}_{(i,j)}+k_m)/2\le
ij/2+Cn'$ then $k_m\le ij/2$, and in the same way $k_M\ge ij/2$.
\end{proof}

This completes the proof of Proposition~\ref{retbound}.


\begin{thebibliography}{9}
\bibitem{Be} M.~Bestvina. Questions in geometric group theory.
Available at \texttt{http://www.math.utah.edu/\~{}bestvina}.
\bibitem{Bo} O.V.~Bogopol'ski\u{i}. Infinite commensurable hyperbolic
groups are bi-Lip\-schitz equivalent. \emph{Algebra and Logic},
36(3):155--63, 1997.
\bibitem{CR} S.~Cleary and T.R.~Riley. A finitely presented group with
unbounded dead-end depth. \emph{Proceedings of the American
Mathematical Society}, 134(2):343--9, 2006.
\texttt{arXiv:math.GR/0406443}.
\bibitem{CT} S.~Cleary and J.~Taback. Dead end words in lamplighter
groups and other wreath products. \emph{The Quarterly Journal of
Mathematics}, 56(2):165--78, 2005. \texttt{arXiv:math.GR/0309344}.
\bibitem{L} J.~Lehnert. Some remarks on depth of dead ends in groups.
Preprint. \texttt{arXiv:math.GR/0703636}, 2007.
\bibitem{RW} T.R.~Riley and A.D.~Warshall. The unbounded dead-end depth
property is not a group invariant. \emph{International Journal of
Algebra and Computation}, 16(5):969--83, 2006.
\texttt{arXiv:math.GR/0504121}.
\bibitem{W} A.D.~Warshall. Deep pockets in lattices and other groups.
Preprint. \texttt{arXiv:math.GR/0611575}, 2007.
\end{thebibliography}
\end{document}